%% file: FCAA_850_Sturm-Liouville-ed.tex
   \documentclass[twoside,reqno,11pt]{fcaa-var} %
   \input fcaa_style_DG %
   \usepackage{hyperref}
  \def\theequation{\arabic{section}.\arabic{equation}}

 \def\themycopyrightfootnote{\vspace*{3pt}
 \copyright \, Year\,  Diogenes  Co., Sofia
   \par  \noindent pp. xxx--xxx, DOI: ......................
   \hfill  \vspace*{-36pt}
  }

\usepackage{mathtools}
\usepackage{amsfonts}

\newcommand{\Real}{\mathbb R}
\newcommand{\Res}{\mathrm{Res}}
\renewcommand{\Re}{\mathrm{Re}}
\renewcommand{\Im}{\mathrm{Im}}

\def\Xint#1{\mathchoice
   {\XXint\displaystyle\textstyle{#1}}%
   {\XXint\textstyle\scriptstyle{#1}}%
   {\XXint\scriptstyle\scriptscriptstyle{#1}}%
   {\XXint\scriptscriptstyle\scriptscriptstyle{#1}}%
   \!\int}
\def\XXint#1#2#3{{\setbox0=\hbox{$#1{#2#3}{\int}$}
     \vcenter{\hbox{$#2#3$}}\kern-.5\wd0}}

\def\dashint{\Xint-}

 \title[SHARP ASYMPTOTICS IN A FRACTIONAL \dots] 
 {SHARP ASYMPTOTICS IN A FRACTIONAL \\ [3pt] STURM-LIOUVILLE PROBLEM }
 \author[\normalsize P. Chigansky, M. Kleptsyna]{\normalsize Pavel Chigansky $^1$, Marina Kleptsyna $^2$}

\begin{document}

 \vbox to 1.5cm { \vfill }

 \bigskip \medskip

 \begin{abstract}

The current research of fractional Sturm-Liouville boundary value problems
focuses on the qualitative theory and numerical methods, and much progress
has been recently achieved in both directions. The objective of this
paper is to explore a different route, namely, construction of explicit
asymptotic approximations for the solutions. As a study case, we consider
a problem with left and right Riemann-Liouville derivatives, for which our
analysis yields asymptotically sharp estimates for the sequence of eigenvalues and
eigenfunctions.

 \medskip


 \smallskip


 \end{abstract}

 \maketitle

 \vspace*{-24pt}


\section{Introduction}

\setcounter{section}{1}
\setcounter{equation}{0}\setcounter{theorem}{0}

The theory of boundary value problems of Sturm-Liouville type is a current subject of research in
fractional calculus and its applications. A number of formulations have been considered, corresponding
to various nonequivalent types of fractional derivatives and motivated by different applications.
Beyond a few special cases, such as e.g. \cite{ZK13}, these problems do not have explicit solutions,
which leaves much space for qualitative theory, numerical methods and asymptotic analysis.

The research in qualitative theory is concerned with the properties of solutions in problems with different
fractional operators and boundary conditions, \cite{KA13},  \cite{K19}, \cite{DA2019}, \cite{LQ15}, \cite{KB15}.
They can be studied by a variety of tools, including functional analysis methods \cite{LQ19}, perturbation techniques \cite{QiChen2011}, variational characterizations \cite{KOM14}, integral transforms \cite{DM2020}. Numerical analysis of fractional Sturm-Liouville problems can also be approached in a number of ways, \cite{BCLK11}, \cite{BC14}, \cite{BC16}, \cite{KCB18}, \cite{OBB2020}.
These references are only a small part of works on the subject in the recent years.

In this paper we address the problem of asymptotic approximation of solutions to fractional eigenvalue problems.
This is a classical theme in functional analysis and mathematical physics, which, in
the fractional setting, has not yet been explored in depth. Such approximations have numerous applications, as they determine properties of various physical quantities, see e.g. \cite{Gibbs69}, \cite{LiS01}, \cite{Br03a}.
Moreover, they can be a decent alternative to numerical methods, which typically remain efficient only for a few dozens of
the first eigenvalues and eigenfunctions. Our computational experiments indicate that even in this range asymptotic approximations may provide a very reasonable accuracy, see Figures 1-2 below.

Our main objective is to draw attention to a technique, which can be useful for deriving
sharp asymptotic estimates of the eigenvalues and eigenfunctions in fractional
Sturm-Liouville problems on a finite interval.
As in the analysis of the fractional ODEs on the whole semi-axis, see e.g. \cite{KST06}, it is based on the
Laplace transform, however in a completely different way. The crucial property which makes the Laplace transform
work on unbounded domains is that it converts the action of fractional operators into multiplication by power
functions. This feature is no longer available, when the domain of functions in question
is bounded, and instead, the method of this paper makes use of its analytic structure instead.

Our approach is capable of producing uniform approximations for various spectral problems, but the detailed analysis
and the ultimate results are quite sensitive to peculiarities of the concrete problem at hand.
As our study case, we will consider a basic, yet nontrivial fractional boundary value problem
\begin{equation}\label{P}\tag{P}
\begin{aligned}
&
 D^\alpha_{1-} D^\alpha_{0+} f (x) = \lambda f(x), \quad x\in [0,1],
\\
&
f(0)=f(1)=0,
\end{aligned}
\end{equation}
with the left and right Riemann-Liouville derivatives of order $\alpha\in (\frac 1 2,1)$,
\begin{align*}
 D^\alpha_{0+}f (x) = & \frac{1}{\Gamma(1-\alpha)}\frac d{dx} \int_0^x (x-t)^{-\alpha}f(t)dt, \\
 D^\alpha_{1-}f (x) = & \frac{1}{\Gamma(1-\alpha)}\frac d{dx} \int_x^1 (t-x)^{-\alpha}f(t)dt.
\end{align*}
These and other standard definitions and formulas from  fractional calculus
can be found in, e.g., \cite{KST06}.

Problem \eqref{P} was studied in \cite{JL16} within the functional analytic framework.
The authors prove that the inverse operator is self-adjoint and compact in suitable spaces.
Consequently its spectrum is discrete and thus the problem has countably many solutions
$(\lambda_n, f_n)_{n\in \mathbb N}$. The eigenvalues $\lambda_n$
are real and nonnegative and, accumulate at infinity, and the corresponding eigenfunctions $f_n$ are
continuous and form an orthogonal basis in the relevant Hilbert space.
Alternatively these results can be obtained by means of reduction to the spectral problem for a
compact self adjoint integral operator with a certain symmetric continuous kernel, see  \cite{KB14}, \cite{KB15}.

\section{Main result} 

\setcounter{section}{2}
\setcounter{equation}{0}\setcounter{theorem}{0}

Our main result is the following theorem, which details asymptotic structure of solutions to the eigenproblem
formulated in the previous section.

\vspace*{-4pt}

\begin{theorem}\label{thm-main}\,
Let $(\lambda_n,f_n)_{n\in \mathbb N}$ be the solutions to \eqref{P} with $\alpha\in (\frac 1 2,1)$, ordered
so that the eigenvalues $\lambda_n$ form a nondecreasing sequence.

\medskip
\noindent
(a) The sequence of frequencies $\rho_n := \lambda_n^{1/(2\alpha)}$ has the asymptotics
\begin{equation}\label{rho_n}
\rho_n =  \pi n + \frac \pi 2 \Big(1-\frac 1 { \alpha}\Big)+O(n^{-1}), \quad n\to\infty.
\end{equation}

\medskip
\noindent
(b) The corresponding eigenfunctions with the unit $L_2$-norm satisfy
\begin{align}\label{fn}
& f_n (x) =
 \sqrt 2\sin \left(  \rho_n x+\frac \pi 4 (1-\alpha)\right) \\
 &\nonumber
+ \int_0^\infty  \Upsilon_0(t)e^{- \rho_n tx} dt
+  (-1)^n\! \int_0^\infty  \Upsilon_1(t)  e^{-  \rho_n t(1-x)} dt
 + r_n(x) n^{-1},\ \, d x\in [0,1],
\end{align}
where  $\sup_n \|r_n\|_\infty<\infty$ and
$\Upsilon_0(t)$ and $\Upsilon_1(t)$ are explicit functions (see \eqref{Upsilons} below).
\end{theorem}

\smallskip  

Several comments are in order.

\medskip

a) Asymptotic approximation of the eigenvalues, implied by this result, is exact up to the second term
with a sharp estimate for the residual,
$$
\lambda_n = (\pi n)^{2\alpha} + \pi (\alpha-1 )(\pi n)^{2\alpha-1} + O(n^{2\alpha-2}), \quad n\to\infty.
$$
The first term of this asymptotics can also be derived through reduction to the eigenproblem  for
an integral operator  \cite{KB14}  and the general results in spectral theory \cite{BS70}.
Our approach does not appeal to such a reduction and is based on direct analysis of the fractional operators.

Formula \eqref{fn} implies that, for large $n$, the eigenfunctions behave as pure harmonics, away from the boundary points
$\{0,1\}$. The two integral terms form the boundary layer, as their contribution is negligible in the
interior of the interval and they force the eigenfunctions to vanish at its endpoints.
It should be emphasised that the boundary layer terms are asymptotically negligible with respect to approximation
in the $L_2([0,1])$ norm, but not in the uniform norm. The second order term in the frequencies $\rho_n$
is non-negligible for the purpose of approximating the eigenfunctions in $L_2([0,1])$ norm.

For $\alpha=1$ problem  \eqref{P} reduces to the classical problem
\vskip -12pt
\begin{align*}
 -  &\frac {d^2}{dx^2}f(x)=\lambda f(x), \\[2pt]
& f(0)=f(1)=0,
\end{align*}
for which elementary calculations give
\vskip -10pt
\begin{equation}\label{classics}
\lambda_n = (\pi n)^2\quad \text{and}\quad f_n(x) = \sqrt 2 \sin(\pi n x), \quad n =1,2,...
\end{equation}
It can be seen that $\Upsilon_0(\cdot)$ and $\Upsilon_1(\cdot)$ in \eqref{fn} vanish as $\alpha\to 1$ and
hence these formulas coincide  formally with \eqref{classics}.

\medskip

b) Numerical approximations to solutions of fractional Sturm-Liouville problems, even as basic as problem \eqref{P},
typically provide good accuracy for a few first eigenvalues and eigenfunctions, but become unstable
already for  $n\ge 30$. The approximation provided by Theorem \ref{thm-main} can thus be a
reasonable alternative. Our numerical experiments indicate that formulas \eqref{rho_n}-\eqref{fn} with
residuals being truncated, turn out to be quite accurate, at least beyond several first eigenpairs,
see Figures 1-2. They also demonstrate significant improvement due to the second order term in
\eqref{rho_n}.

\begin{center}
\includegraphics[scale=0.4]{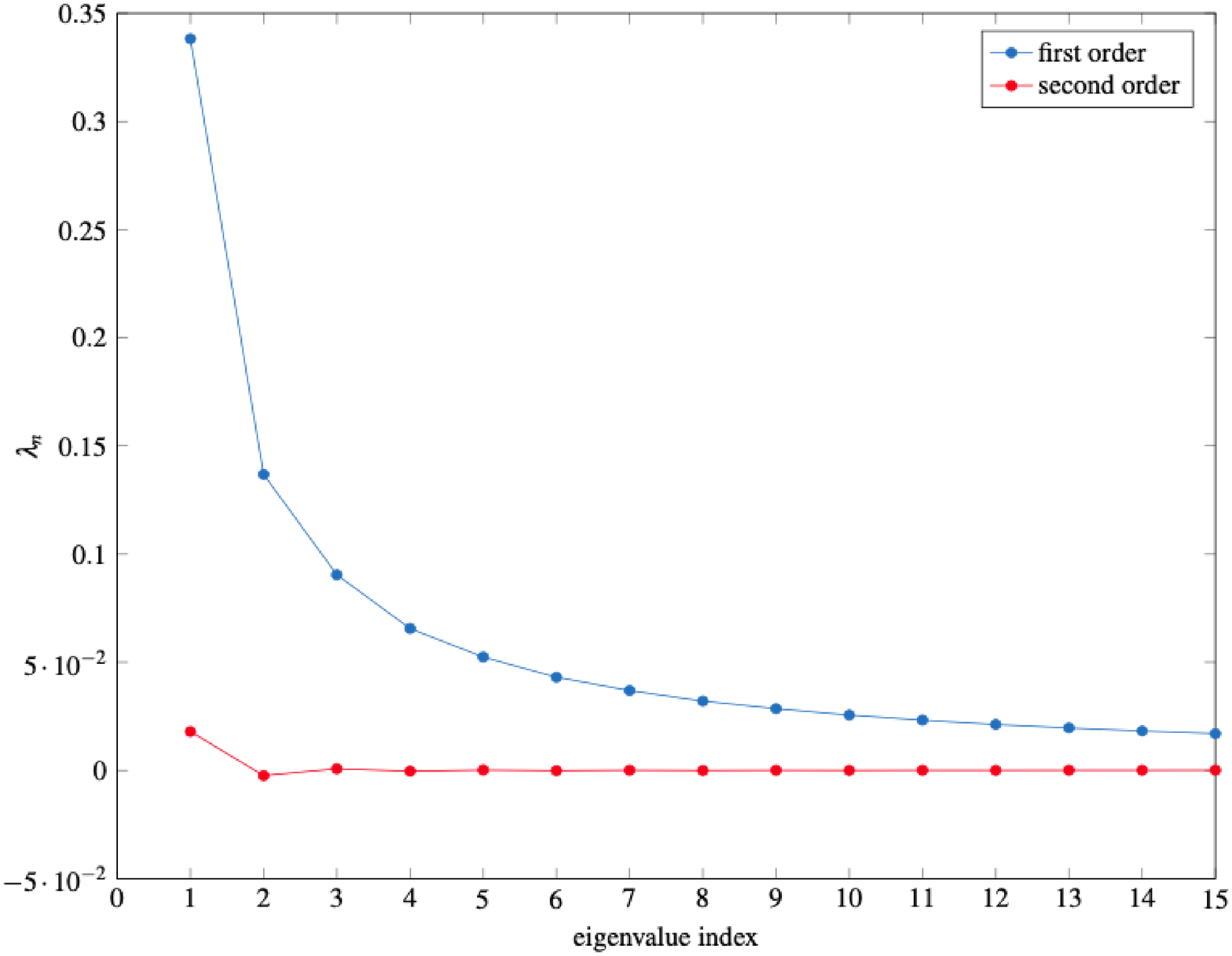}

\smallskip 
 \begin{minipage}{0.8\textwidth} Figure 1. The relative error $\widehat \lambda_n/\widetilde \lambda_n^{(i)}-1$ versus $n$
for $\alpha= 3/4$, where $\widehat \lambda_n$ is a high precision numerical approximation of the
eigenvalues $\lambda_n$ and $$\widetilde \lambda_n^{(1)} = (\pi n)^{2\alpha}\quad \text{and} \quad
\widetilde \lambda_n^{(2)} = \big(\pi n+\tfrac \pi 2  (1-\tfrac 1 { \alpha} )\big)^{2\alpha}$$
are the  first and second order  asymptotic  approximations.
\end{minipage}
  \end{center}

\begin{center}
  \includegraphics[scale=0.4]{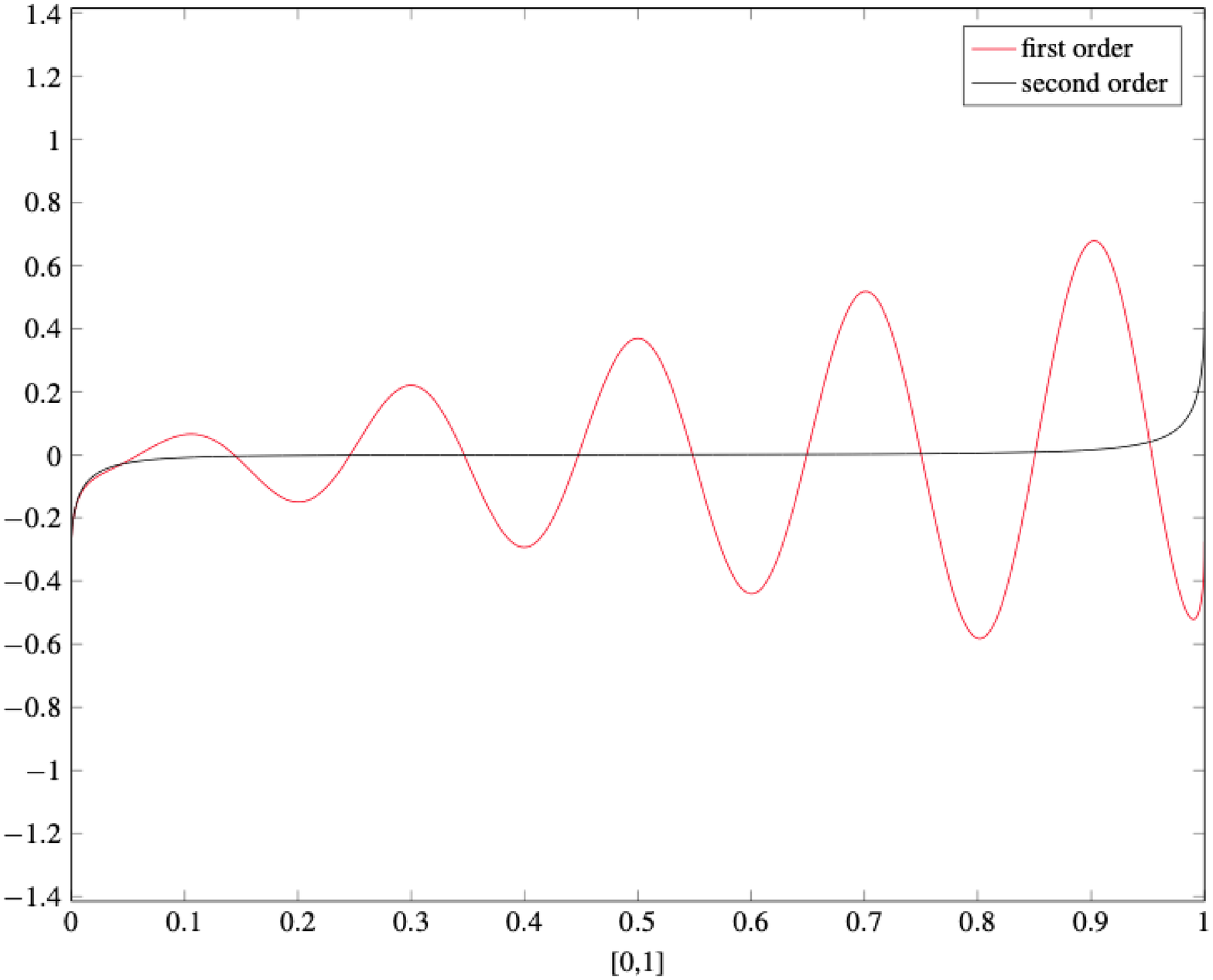}

 \smallskip 
   \begin{minipage}{0.8\textwidth} Figure 2. The error $\widehat f_n(x) -\widetilde f_n^{(i)}(x)$ versus $x\in [0,1]$ for  $\alpha=3/4$
and $n=10$, where $\widehat f_n$ is a high precision numerical approximation of the eigenfunction $f_n$ and
$$
\widetilde f_n^{(i)}  =  \sqrt 2\sin \big( \widehat \rho_n^{(i)} x+\tfrac \pi 4 (1-\alpha)\big), \quad i=1,2,
$$
with $\widehat\rho_n^{(i)} = (\widehat \lambda^{(i)})^{1/(2\alpha)}$, are the first and second order approximations
for the oscillating part of the eigenfunctions. The effect of truncated boundary layer terms is clearly visible for the
second order approximation.
\end{minipage}
  \end{center}

\medskip

c) Problem \eqref{P} is closely related to a slightly different problem
\vskip -10pt
\begin{equation}\label{C}\tag{P'}
\begin{aligned}
&
 \prescript{c}{}D^\alpha_{1-} \prescript{c}{}D^\alpha_{0+} u (x) = \lambda u(x), \quad x\in [0,1],
\\[2pt] 
&
u(0)=0,\ \prescript{c}{}D^\alpha_{0+} u(1)=0,
\end{aligned}
\end{equation}
with the left and right Caputo derivatives
\begin{align*}
\prescript{c}{}D^\alpha_{0+} f(x) :=\, & D^\alpha_{0+} \big(f(x)-f(0)\big),\\
\prescript{c}{}D^\alpha_{1-} f(x) :=\, & D^\alpha_{1-} \big(f(x)-f(1)\big).
\end{align*}
Reduction to the integral operators as in \cite{KB14}, \cite{KB15} reveals that \eqref{C} with $\alpha \in (0,1)$
is equivalent to the problem
$$
\int_0^1 K(x,y)f(y)dy =  \lambda^{-1} f(x), \quad x\in [0,1],
$$
with the kernel
\vskip -10pt
\begin{equation}\label{Kxy}
K(x,y) =\frac 1 {\Gamma(\alpha)^2} \int_0^{x\wedge y} (x-y)^{\alpha-1}(y-t)^{\alpha-1}dt,
\end{equation}
while problem \eqref{P} is equivalent to
\begin{equation}\label{PI2}
\int_0^1 \Big(K(x,y)-\frac{K(x,1)K(1,y)}{K(1,1)}\Big)f(y)dy =  \lambda^{-1} f(x), \quad x\in [0,1],
\end{equation}
with $K(x,y)$ as in \eqref{Kxy}.

Our approach applies to  problem \eqref{C} with minor adjustments and it can be seen that the assertion of
Theorem \ref{thm-main} remains valid, this time for $\alpha \in (0,1)$, with
\vskip -10pt 
$$
\rho_n = \pi n - \frac \pi 2 + O(n^{-1}),
$$
and different but still explicit functions $\Upsilon_j(t)$. Curiously, in this case
the solutions can be shown to have very particular value at the right endpoint,
$$
\big|f_n(1)\big| =  \sqrt{2\alpha},
$$
at least for all sufficiently large $n$.

In the theory of stochastic processes,  kernel \eqref{Kxy} with $\alpha \in (\frac 1 2,1)$ is the covariance
function of the so called Riemann-Liouville process, one of the two most common fractional generalisations of
the standard Brownian motion. The kernel in \eqref{PI2} with $\alpha \in (\frac 1 2, 1)$ is the covariance
function of the corresponding bridge,  obtained by conditioning the Riemann-Liouville process to vanish at $x=1$.

For $\alpha \in (0,\frac 1 2]$ the covariance function \eqref{Kxy} is still well defined, but it remains continuous
only off diagonal, where it has a weak singularity. Such covariance function does not correspond to a stochastic
process in the usual sense, as it cannot have continuous paths. Instead it can be interpreted as covariance
of the fractional noise process, a formal derivative of the Riemann-Liouville process. This can be made precise by
e.g., considering stochastic integrals, see \cite{ChK}.  Obviously, the bridge kernel as in \eqref{PI2} can no longer
be defined in this case. In fact, our analysis implies that problem \eqref{P} with $\alpha \in (0,\frac 1 2)$
cannot have more than finitely many continuous solutions.

\section{Proof of Theorem \ref{thm-main}}\label{sec-proof}

\setcounter{section}{3}
\setcounter{equation}{0}\setcounter{theorem}{0}

The proof is based on the approach to analysis of integral equations with weakly
singular kernels, introduced in \cite{Ukai}.
In its original form this method applies to operators with difference kernels, see also \cite{P74}, \cite{P03}.
Recently it was extended to various forms of {\em integrated} difference kernels, which allowed to compute the exact
spectral asymptotics of covariance operators for a number of stochastic processes, including the fractional
Brownian motion, \cite{ChK}, \cite{ChKM20}, \cite{ChKM-AiT}, \cite{N19}.

While inversion of the fractional derivatives in problem \eqref{P} still yields an integral operator, its kernel
does not seem to have any simple difference structure.  The principal contribution of this paper is generalisation of
the method to the integro-differential operators of the fractional types, which makes it potentially applicable to
corresponding Sturm-Liouville problems.
To avoid repetitions and focus only on the new elements of the proof, we will omit calculations, which can be found
elsewhere and provide with exact references.

\subsection{Conventions and notations}

The proof uses some standard notations and classical tools from complex analysis.
Unless stated otherwise, we will use the conventional principal branches for multivalued functions.
Often we will encounter functions, which are holomorphic on the complex plane, cut along the real
(semi-) axis with a finite jump discontinuity across the cut.
For such a {\em sectionally holomorphic} function $\Psi(z)$, we will use the limit notations
$$
\Psi^+(t) = \lim_{\Im(z)>0, z\to t} \Psi(z)\quad \text{and}\quad \Psi^-(t) = \lim_{\Im(z)<0, z\to t} \Psi(z),
\quad t\in \Real.
$$

Finding a function $\Psi(z)$, sectionally holomorphic on $\mathbb{C}\setminus \Real_+$ and satisfying the boundary condition
$$
\Psi^+(t)-\Psi^-(t) = g(t), \quad t\in \Real_+,
$$
where $g(\cdot)$ is a given H\"older function on $\Real_+\cup\{\infty\}$, is known as the Hilbert boundary value
problem. The solution is given by the Sokhotski-Plemelj formula
$$
\Psi(z) = \frac 1 {2\pi i} \int_0^\infty \frac{g(t)}{t-z}dt + P(z), \quad z\in \mathbb{C}\setminus\Real_+,
$$
where $P(\cdot)$ is a polynomial, which matches the growth of $\Psi(\cdot)$ at infinity.
A comprehensive account of the related theory can be found in \cite{Gahov}.

\subsection{A brief preview} 

In a nutshell the main idea of the method is to reduce  eigenproblem \eqref{P} to an
equivalent integro-algebraic system of equations, more tractable for asymptotic analysis.
In Lemma \ref{lem-L} below we show that the Laplace transform
\vskip -10pt
\begin{equation}\label{Laplace}
\widehat f(z) = \int_0^1 e^{-zx}f(x)dx, \quad z\in \mathbb{C},
\end{equation}
 of a solution $f$ to problem \eqref{P}
satisfies the representation
\begin{equation}\label{repre}
\widehat f(z)= - \frac { c_f}\lambda e^{-z}
-  \frac {\Phi_0(z)+e^{-z}\Phi_1(-z) }
{\Lambda(z)}, \quad z\in \mathbb{C},
\end{equation}
where  functions $\Phi_0(z)$ and $\Phi_1(z)$ are  sectionally holomorphic on $\mathbb C\setminus \Real_+$
and $c_f$ is a functional of $f$, constant in $z$.
The function $\Lambda(z)$ is defined by a closed form expression, which vanishes only at a pair of purely imaginary conjugate
zeros and is sectionally holomorphic on the cut plane $\mathbb C\setminus \Real$.

Since a priori the Laplace transform $\widehat f(z)$ is an entire function, all singularities in \eqref{repre}
must be removable. Removing discontinuity along the real line yields   representations of
$\Phi_0(z)$ and $\Phi_1(z)$ in terms of solutions to an auxiliary system of integral equations on $\Real_+$,
see Lemma \ref{lem:34}. Removal of the poles produces an algebraic condition, which binds together values of these
functions at certain points in the complex plane. This integro-algebraic system of equations, described in
Lemma \ref{lem35}, is shoen to have countably many solutions, whose asymptotic analysis leads to the expressions
claimed in Theorem \ref{thm-main}.

Thereby problem \eqref{P} reduces to finding a pair of functions, sectionally holomorphic on $\mathbb C\setminus\Real_+$
with a given jump discontinuity on $\Real_+$ and specific growth near the origin and at infinity.
Remarkably, this seemingly more complicated problem turns out to be amenable to asymptotic analysis.

\subsection{The Laplace transform} 

Our starting point is the following representation for the Laplace transform.
\vspace*{-3pt}

\begin{lemma}\label{lem-L}
Let $(\lambda, f)$ be a solution to \eqref{P}. The Laplace transform \eqref{Laplace} of $f$
satisfies the representation
\begin{equation}\label{hatfrep}
\widehat f(z)=- \frac { c_f}\lambda e^{-z}
- \frac {\Phi_0(z)+e^{-z}\Phi_1(-z)}
{\Lambda(z)}, \quad z\in \mathbb C,
\end{equation}
where $c_f =  \lim_{x\to 1}I_{1-}^{1-\alpha} D^\alpha_{0+} f(x)$,
\begin{equation}\label{Lambda}
\Lambda(z) = \frac \pi {\sin (\alpha \pi )} \big( z^{\alpha } -  \lambda   (-z)^{-\alpha}\big),
\end{equation}
and functions $\Phi_0(z)$ and $\Phi_1(z)$, defined in \eqref{psiz} below,
are sectionally holomorphic on $\mathbb C\setminus \Real_+$ and satisfy the growth estimates
\begin{equation}\label{Phi_est_inf}
\Phi_0(z) =O(z^{-\alpha})\quad \text{and}\quad \Phi_1(z)=O(z^\alpha), \quad z\to \infty,
\end{equation}
\vskip -3pt \noindent
 and
 \vskip -11pt
\begin{equation}\label{Phi_est_0}
\Phi_0(z) = O(z^{- \alpha })\quad \text{and}\quad \Phi_1(z)=O(1), \quad z\to 0.
\end{equation}
\end{lemma} 

\proof
Applying the left Riemann-Liouville integral
\vskip -10pt
$$
\big(I^\alpha_{1-} f\big)(x)= \frac 1 {\Gamma (\alpha)}  \int_x^1(t-x)^{\alpha-1}f(t)dt
$$
to both sides of the equation in \eqref{P} and using the composition rule \cite[Lemma 2.6]{KST06}, we can write
\vskip -10pt
$$
 D^\alpha_{0+} f (x) - \frac {c_f}{\Gamma(\alpha)}  (1-x)^{\alpha-1}= \lambda  I^\alpha_{1-} f (x).
$$
Following \cite{N19}, define the function
\vskip -12pt
\begin{equation}\label{psi}
\psi(x) = \int_x^1 f(y)dy +\frac {c_f}\lambda,
\end{equation}
\vskip -4pt \noindent 
then
\begin{align*}
&
\lambda  I^\alpha_{1-} f (x) + \frac {c_f}{\Gamma(\alpha)} (1-x)^{\alpha-1} 
\\
&
= -\frac{\lambda }{\Gamma(\alpha)}\int_x^1 \psi'(y)(y-x)^{\alpha-1}dy
+ \frac {c_f}{\Gamma(\alpha)} (1-x)^{\alpha-1} \\
&
=  \frac{\lambda }{\Gamma(\alpha)}\frac 1 \alpha \frac d{dx}\int_x^1 \psi'(y)(y-x)^{\alpha}dy
+ \frac {c_f}{\Gamma(\alpha)} (1-x)^{\alpha-1} \\
&
= \frac{\lambda \psi(1)}{\Gamma(\alpha)\alpha}    \frac d{dx}
 (1-x)^{\alpha}  - \frac{\lambda }{\Gamma(\alpha)}  \frac d{dx} \int_x^1 \psi(y)(y-x)^{\alpha-1}dy
+ \frac {c_f}{\Gamma(\alpha)} (1-x)^{\alpha-1} \\
&
=  - \frac{\lambda }{\Gamma(\alpha)}  \frac d{dx} \int_x^1 \psi(y)(y-x)^{\alpha-1}dy.
\end{align*}
Thus  problem \eqref{P} takes the equivalent form
\begin{equation}\label{prpr}
\begin{aligned}
&
\frac 1 {\Gamma(1-\alpha)}\frac d{dx} \int_0^x \psi'(y) (x-y)^{-\alpha}dy
 =    \frac{\lambda }{\Gamma(\alpha)}  \frac d{dx} \int_x^1 \psi(y)(y-x)^{\alpha-1}dy,  \\[3pt]
&
\psi'(0)=\psi'(1)=0.
\end{aligned}
\end{equation}
Using the identity
\vskip -10pt
$$
s^{-r} = \frac 1 {\Gamma(r)} \int_0^\infty t^{r-1} e^{-st}dt, \quad r >0, \ s>0,
$$
the equation in \eqref{prpr} can be written as
\vskip -10pt
\begin{equation}\label{v0u0}
 v_0'(x)=  \lambda   u_0'(x),
\end{equation}
where we denoted
\vskip -10pt
$$
u_0 (x) =  \int_0^\infty t^{-\alpha} u(x,t) dt,\quad
u(x, t) = \int_x^1 \psi(y)  e^{-t(y-x) }dy,
$$
and
$$
v_0 (x) = \int_0^\infty t^{\alpha-1} v(x,t) dt,\quad
v(x,t) = \int_0^x \psi'(y) e^{-t(x-y)}dy.
$$
Differentiating with respect to $x$ shows that $u(x,t)$ solves the differential equation
\vskip -12pt
\begin{align*}
u'(x,t) & =  tu(x,t)- \psi(x),   \\
u(1,t) & = 0.
\end{align*}
Applying the Laplace transform we get
$$
\widehat u'(z,t) = t \widehat u(z,t) -\widehat \psi(z).
$$
On the other hand, $\widehat u'(z,t) =  -u(0,t)+z\widehat u(z,t)$, and therefore
$$
\widehat u(z,t) = \frac{\widehat \psi(z)-u(0,t)}{t-z}.
$$
\vskip -3pt \noindent
It follows that
\vskip -12pt
\begin{equation}\label{hatu0}
\widehat u_0(z) =
\widehat \psi(z) M_u(z)
-  \int_0^\infty  \frac{ t^{-\alpha} }{t-z}u(0,t)  dt,
\end{equation}
\vskip -3pt \noindent
where  
\vskip -10pt
$$
M_u(z):=\int_0^\infty  \frac{t^{-\alpha}  }{t-z}  dt =\frac \pi {\sin (\alpha \pi)}  (-z)^{-\alpha}.
$$
The last equality holds due to the identity
$$
\int_0^\infty \frac{t^{-\beta}}{t+z}dt=\frac \pi {\sin (\pi \beta)}  z^{-\beta}, \quad \beta \in (0,1), \ \arg(z)\in [-\pi,\pi).
$$
Similarly,
\vskip -12pt
\begin{align*}
v'(x,t) = &  -tv(x,t)+\psi'(x),  \\
v(0,t)= & 0,
\end{align*}
and, since $\widehat v'(z,t) = v(1,t)e^{-z} + z\widehat v(z,t)$, we get
$$
\widehat v(z,t) =   \frac{\widehat \psi'(z) -v(1,t)e^{-z}}{t+z}=\frac{ \psi(1)-v(1,t) }{t+z}e^{-z}
+ \frac{  z\widehat \psi(z) -\psi(0)  }{t+z},
$$
and consequently,
\vskip -11pt
\begin{equation}
\begin{aligned}
\label{hatv0}
\widehat v_0(z)=\, & e^{-z}\psi(1) M_v(-z)
 -  e^{-z}\int_0^\infty    \frac{t^{\alpha-1} }{t+z} v(1,t)
 dt \\
&
+ \big(z\widehat \psi(z) -\psi(0) \big)M_v(-z),
\end{aligned}
\end{equation}
where we defined
\vskip -10pt
$$
M_v(z)=\int_0^\infty  \frac{   t^{\alpha-1} }{t-z} dt = \frac \pi {\sin (\alpha \pi )}  (-z)^{\alpha-1}.
$$

In the Laplace domain, the equation \eqref{v0u0} is equivalent to
$$
 z\big(\widehat v_0(z) -\lambda  \widehat u_0(z) \big)
=
c\big(1-e^{-z}\big),
$$
where $c:=  v_0(1) =  -\lambda u_0(0)$. Plugging \eqref{hatu0} and \eqref{hatv0} into this equation
and rearranging we arrive at
\vskip -12pt 
\begin{equation}\label{psiz}
 z\widehat \psi(z) -\psi(0)
= \frac { \Phi_0(z)+e^{-z}\Phi_1(-z) } {\Lambda(z)},
\end{equation}
where  
\vskip -10pt
\begin{equation}\label{Phi01}
\begin{aligned}
\Phi_0(z)  =\,
&   \psi(0)\lambda M_u(z)-\lambda
\int_0^\infty    \frac{ t^{1-\alpha} }{t-z} u(0,t) dt,
 \\
\Phi_1(z)  =
&
 \psi(1) z  M_v(z) -  \int_0^\infty    \frac{t^{\alpha } }{t-z} v(1,t) dt,
\end{aligned}
\end{equation}
\vskip -4pt \noindent 
and
\vskip -10pt
$$
\Lambda(z) = zM_v(-z)-  \lambda M_u(z) =
\frac \pi {\sin (\alpha \pi )} \big( z^{\alpha } -  \lambda   (-z)^{-\alpha}\big).
$$
Representation \eqref{hatfrep} follows from \eqref{psiz} since $z\widehat \psi (z)-\psi(0)=-\widehat f(z)-e^{-z} \psi(1)$.
The growth estimates \eqref{Phi_est_0} and \eqref{Phi_est_inf} are verified by standard calculations.
\proofend

\medskip

The function $\Lambda(z)$ largely determines the structure of the problem.
The following lemma summarises some of its key properties, relevant to further analysis.

\vspace*{-3pt}

\begin{lemma}\label{lem-Lambda}\

\medskip
\noindent
a)
The function $\Lambda(z)$  in \eqref{Lambda} is sectionally holomorphic on $\mathbb C\setminus \Real$ and
discontinuous across the real line with the limits
$$
\Lambda^\pm(t) =
\frac \pi {\sin (\alpha \pi )}
\begin{cases}
t^{\alpha } -  \lambda   e^{\pm \alpha \pi  i}t^{-\alpha}, & t> 0, \\
e^{\pm \alpha \pi i }|t|^{\alpha } -  \lambda  |t|^{-\alpha}, & t<0,
\end{cases}
$$
which satisfy the symmetries
\begin{equation}\label{Lambdapm}
\begin{aligned}
& \Lambda^+(t) = \overline{\Lambda^-(t)}, \\
&
\Lambda^+(t) = e^{\alpha\pi i}\Lambda^-(-t),\\
&
\frac {\Lambda^+(t)}{\Lambda^-(t)} = e^{2\alpha \pi i}\frac{\Lambda^-(-t)}{\Lambda^+(-t)}.
\end{aligned}
\end{equation}

\medskip
\noindent
b)
As $t$ varies from $0$ to $\infty$, the angle $\theta(t)=\arg\{\Lambda^+(t)\}$ increases continuously from
$\theta(0+) = (\alpha-1)\pi <0$ to $\theta(\infty)=0$ and
$$
|\theta(t)|=O(t^{-2\alpha}), \quad \text{as}\ \to \infty.
$$

\medskip
\noindent
c) $\Lambda(z)$ vanishes only at simple zeros $\pm z_0=\pm i \rho$ with
$$
\rho = \lambda^{1/(2\alpha)},
$$
and the function $\theta_0(t) := \theta(\rho t)$ does not depend on $\rho$.
\end{lemma} 

\proof
All the claims are verified by direct calculations.
\proofend

\subsection{Removal of singularities} 

Since the integration in \eqref{Laplace} is over a finite interval, $\widehat f(z)$ is an entire
function and therefore all singularities in \eqref{hatfrep} must be removable.
The discontinuity is removed by equating the limits  in the upper and lower half planes,
\begin{align*}
&
\frac 1{\Lambda^+(t)}\Big(e^{-t}\Phi_1(-t)+\Phi_0^+(t)\Big) = \frac 1 {\Lambda^-(t)}\Big(e^{-t}\Phi_1(-t)+ \Phi_0^-(t)\Big), \quad t\in \Real_+,
\\
&
\frac 1{\Lambda^+(t)}\Big(e^{-t}\Phi_1^-(-t)+\Phi_0(t)\Big) = \frac 1 {\Lambda^-(t)}\Big(e^{-t}\Phi_1^+(-t)+ \Phi_0(t)\Big), \quad t\in \Real_-,
\end{align*}
which, in view of \eqref{Lambdapm}, can be written as
\begin{equation}\label{bndcnd}
\begin{aligned}
&
\Phi_0^+(t)-  \frac {\Lambda^+(t)} {\Lambda^-(t)}\Phi_0^-(t)  = e^{-t}\Phi_1(-t)\Big(\frac {\Lambda^+(t)} {\Lambda^-(t)} - 1\Big),
\\
&
   \Phi_1^+(t)-e^{-2\alpha \pi i}\frac{\Lambda^+(t)}{\Lambda^-(t)}   \Phi_1^-(t) = e^{-t}\Phi_0(-t)
   \Big(1-e^{-2\alpha \pi i}\frac{\Lambda^+(t)}{\Lambda^-(t)}  \Big),
\end{aligned}
\quad t\in \Real_+.
\end{equation}
\vskip -3pt \noindent 
Removal of the poles in \eqref{hatfrep} implies
\begin{equation}\label{cndpoles}
\Phi_0(z_0)+\Phi_1(-z_0)   e^{-z_0} =0,
\end{equation}
since the zeros of $\Lambda(z)$ are purely imaginary and conjugate,  and $\overline{\Phi_j(z)}=\Phi_j(\overline z)$ by definition \eqref{Phi01}.

\subsection{Auxiliary integro-algebraic system} 

The next step is to show that any pair of functions $\Phi_0(z)$ and $\Phi_1(z)$, which are sectionally holomorphic
on the cut plane $\mathbb C\setminus \Real_+$, satisfy  boundary conditions \eqref{bndcnd} and  growth
estimates \eqref{Phi_est_0}-\eqref{Phi_est_inf}, solve a certain auxiliary system of integral equations. Along with
condition \eqref{cndpoles} these equations form an integro-algebraic system, whose solutions correspond to
solutions of the spectral problem \eqref{P} under consideration.

To this end we will use the classical technique of solving the Hilbert boundary value problems, see e.g. \cite{Gahov}.
Consider the function
\begin{equation}\label{Xcz}
X_c(z) = \exp \left(\frac 1 \pi \int_0^\infty \frac{\theta(t)} {t-z}dt\right), \quad z\in \mathbb{C}\setminus \Real_+,
\end{equation}
where $\theta(t)$ was introduced in Lemma \ref{lem-Lambda}.
By the Sokhotski-Plemelj theorem, $X_c(z)$ is sectionally holomorphic on $\mathbb C\setminus \Real_+$
and its limits satisfy
\begin{equation}\label{XL}
\frac{X_c^+(t)}{X_c^-(t)} = e^{2i \theta(t)} = \frac{\Lambda^+(t)}{\Lambda^-(t)}, \quad t\in \Real_+,
\end{equation}
where the last equality holds by \eqref{Lambdapm}. Some  useful
properties of this function are gathered in the following lemma.

\vspace*{-3pt}

\begin{lemma}\label{lem-Xc0}\
The function  $X_{c0}(z) := X_c(\rho z)$ does not depend on $\rho$, satisfies the growth estimates
\begin{equation}\label{Xcest}
\begin{aligned}
X_{c0}(z)=
\begin{cases}
O(z^{1-\alpha}),  &  \text{as\ } z\to 0,  \\
1 -   b_\alpha z^{-1}  + O(z^{-2\alpha}), & \text{as}\ z\to\infty,
\end{cases}
\end{aligned}
\end{equation}
with
$
b_\alpha = \cot  \left(\frac \pi {2\alpha}\right),
$
and has the following explicit value at $z=i$,
\begin{equation}\label{Xc0val}
X_{c0}(i) = \sqrt{\alpha} \exp\left(-\frac \pi 4 (1-\alpha)i\right).
\end{equation}
\end{lemma} 

\proof
Evaluating \eqref{Xcz} at $\rho z$ and changing integration variable accordingly gives
\vskip -12pt
$$
X_{c0}(z) = \exp \left(\frac 1 \pi \int_0^\infty \frac{\theta_0(t)} {t-z}dt\right),
$$
where $\theta_0(t)$ does not depend on $\rho$, see Lemma \ref{lem-Lambda} (c).
The estimate near the origin in \eqref{Xcest} is obtained by means of integration by parts,
$$
X_{c0}(z) = (-z)^{-\theta_0(0+)/\pi}\exp \left(\frac 1 \pi \int_0^\infty \theta'_0(t)  \log (t-z)dt\right)= O(z^{1-\alpha}), \  z\to 0,
$$
since $\theta_0(0+)=(\alpha-1)\pi$, see Lemma \ref{lem-Lambda} (b).
The estimate at infinity holds  since
\vskip -12pt
\begin{align*}
X_{c0}(z) =\,  & \exp \left(-z^{-1} \frac 1 \pi \int_0^\infty \theta_0(t) dt
+
z^{-1}\frac 1 \pi \int_0^\infty \frac{t\theta_0(t)} {t-z}dt\right) 
\\
&
=  1 - z^{-1}   b_\alpha  + O(z^{-2\alpha}), \quad z\to\infty,
\end{align*}
where 
\vskip -10pt
$$
b_\alpha :=   \frac 1 \pi \int_0^\infty \theta_0(t) d t.
$$
The exact value of the latter integral and the constant in \eqref{Xc0val} are found by direct calculation, which use the explicit expressions for $\Lambda^\pm (t)$ from Lemma \ref{lem-Lambda}, similarly to \cite[Lemma 5.5 ]{ChK}.
\proofend

\medskip 

Define $X(z) :=z^{-1}X_c(z)$ and $Y(z) := (-z)^{\alpha-1} X_c(z)$.
These functions are sectionally holomorphic on $\mathbb C\setminus \Real_+$ and, in view of \eqref{XL}, satisfy
 boundary conditions
 \vskip -10pt
$$
\begin{aligned}
\frac{X^+(t)}{X^-(t)} =\, &  \frac{\Lambda^+(t)}{\Lambda^-(t)}, \\
\frac{Y^+(t)}{Y^-(t)} =\, & e^{-2\alpha\pi i  }\frac{\Lambda^+(t)}{\Lambda^-(t)},
\end{aligned}
\qquad t\in \Real_+.
$$
Plugging these formulas into equations \eqref{bndcnd} and rearranging gives
$$
\begin{aligned}
&
\frac {\Phi_0^+(t)}{X^+(t)}-  \frac{ \Phi_0^-(t) }{X^-(t)} = e^{-t} Y(-t)\Big(\frac{1}{X^-(t)} - \frac 1 {X^+(t)}\Big)\frac{\Phi_1(-t)}{Y(-t)},
\\
&
\frac {\Phi_1^+(t)}{Y^+(t)}-\frac{\Phi_1^-(t)}{Y^-(t)}   = e^{-t}X(-t)
   \Big(\frac 1 {Y^+(t)}-\frac{1}{Y^-(t)} \Big)\frac{\Phi_0(-t)}{X(-t)},
\end{aligned}
\quad t\in \Real_+.
$$
This shows that the functions
\begin{equation}\label{Psi01z}
\Psi_0(z)  = \frac{\Phi_0(z)}{X(z)}\quad \text{and}\quad \Psi_1(z)  = \frac{\Phi_1(z)}{Y(z)},
\end{equation}
are sectionally holomorphic on $\mathbb C\setminus \Real_+$ with the boundary conditions
$$
\begin{aligned}
&
\Psi_0^+(t) -  \Psi_0^-(t) =   2i e^{-t} g(t)\Psi_1(-t),
\\
&
\Psi_1^+(t) -\Psi_1^-(t)   =  2i e^{-t} h(t) \Psi_0(-t),
\end{aligned}
\qquad t\in \Real_+,
$$
where  
\vskip -10pt
\begin{align*}
g(t) := & \frac 1 {2i}   Y(-t)\Big(\frac{1}{X^-(t)} - \frac 1 {X^+(t)}\Big), \\
h(t) := &  \frac 1 {2i} X(-t)\Big(\frac 1 {Y^+(t)}-\frac{1}{Y^-(t)} \Big).
\end{align*}
These two functions turn out to be real valued and satisfy the scaling properties
\vskip-12pt
\begin{align*}
g(\rho t) =\, & \rho^\alpha  t^{\alpha } \sin(\theta_0(t))
\exp \left(-\frac {2t} \pi \dashint_0^\infty \frac{\theta_0(\tau)} {\tau^2-t^2}d\tau  \right) = :\rho^\alpha g_0(t),
\\
h(\rho t) =\, & -
\rho^{-\alpha}  t^{-\alpha } \sin \big(\theta_0(t)-\alpha\pi\big)
\exp \left(-\frac {2t} \pi \dashint_0^\infty \frac{\theta_0(\tau)} {\tau^2-t^2}d\tau  \right) =: \rho^{-\alpha} h_0(t),
\end{align*}
 where the integral is in the sense of Cauchy's principal value and $g_0(\cdot)$ and $h_0(\cdot)$ do not depend on $\rho$.

In view of estimates \eqref{Phi_est_0}  and \eqref{Xcest}, both $\Psi_0(z)$ and $\Psi_1(z)$ are bounded in the vicinity
of zero and therefore, due to the Sokhotski-Plemelj theorem, satisfy
\vskip -10pt 
\begin{equation}\label{PsieqP}
\begin{aligned}
\Psi_0(z) & =   \frac 1 \pi \int_0^\infty \frac{e^{-\tau} g(\tau)}{\tau-z}\Psi_1(-\tau) d\tau + P_0(z),\\
\Psi_1(z) &   =  \frac 1 \pi \int_0^\infty \frac{ e^{-\tau} h(\tau) }{\tau-z}\Psi_0(-\tau)d\tau + P_1(z),
\end{aligned}
\qquad z\in \mathbb C\setminus\Real_+,
\end{equation}
where $P_j(z)$'s are polynomials, whose growth at infinity matches that of $\Psi_j(z)$'s.

Let us estimate the growth of $\Psi_1(z)$ as $z\to\infty$. To this end note that the integral term of $\Phi_1(z)$ in
\eqref{Phi01} satisfies
\begin{align*}
&
z^{1-\alpha}\int_0^\infty    \frac{t^{\alpha } }{t-z} v(1,t) dt=
z^{1-\alpha}|z|^\alpha\int_0^\infty    \frac{t^{\alpha } }{  t-z/|z|} v(1,t|z|) dt  \\
&
=  (z/|z|)^{1-\alpha} \int_0^\infty    \frac{t^{\alpha-1}   }{  t-z/|z|} \int_0^1 \psi'(y) |z| te^{-|z|t(1-y)}dy dt
\xrightarrow{|z|\to\infty}0,
\end{align*}
where the limit holds, since the inner integral converges to $\psi'(1)=0$.
Consequently, as $z\to\infty$,
\begin{align*}
\Psi_1(z)  =\, & \frac{\Phi_1(z)}{Y(z)} =
\psi(1) \frac \pi {\sin (\alpha \pi )}\frac z {  X_c(z)} \big(1+o(1)\big)  \\
&
= \psi(1) \frac \pi {\sin (\alpha \pi )}\big(\rho b_\alpha + z\big) \big(1+o(1)\big), \quad z\to \infty.
\end{align*}
Similar calculations show that $\Psi_0(z)=O(1)$ as $z\to\infty$.
In view of these estimates \eqref{PsieqP} takes the form
$$
\begin{aligned}
\Psi_0(z) & =   \frac 1 \pi \int_0^\infty \frac{e^{-\tau} g(\tau)}{\tau-z}\Psi_1(-\tau) d\tau + c_0,\\
\Psi_1(z) &   =  \frac 1 \pi \int_0^\infty \frac{ e^{-\tau} h(\tau) }{\tau-z}\Psi_0(-\tau)d\tau + c_1 (\rho b_\alpha+z),
\end{aligned}
\qquad z\in \mathbb C\setminus\Real_+,
$$
where $c_0$ and $c_1$ are real constants,
\begin{equation}\label{cf2c1}
c_1 = \psi(1) \frac \pi {\sin (\alpha \pi )} = \frac{c_f}{\lambda} \frac \pi {\sin (\alpha \pi )}.
\end{equation}
In particular, for $z= - t$ we obtain the system of integral equations for $\Psi_0(-t)$ and $\Psi_1(-t)$,
\begin{equation}\label{Psieqt}
\begin{aligned}
\Psi_0(-t) & =   \frac 1 \pi \int_0^\infty \frac{e^{-\tau} g(\tau)}{\tau+t}\Psi_1(-\tau) d\tau + c_0,\\
\Psi_1(-t) &   =  \frac 1 \pi \int_0^\infty \frac{ e^{-\tau} h(\tau) }{\tau+t}\Psi_0(-\tau)d\tau + c_1 (\rho b_\alpha-t),
\end{aligned}
\qquad t\in  \Real_+.
\end{equation}
In view of the growth estimates for $\Psi_0(z)$ and $\Psi_1(z)$, the integrals in the right hand side
define functions in $L_2(\Real_+,\Real)$. Finally,  condition \eqref{cndpoles} reads
\vskip -10pt
\begin{equation}\label{algebr}
\Psi_0(i\rho )X(i\rho)+\Psi_1(-i\rho) Y(-i\rho)  e^{-i\rho} =0.
\end{equation}
\vskip 3pt \noindent
Thus we arrive at the following result.

\vspace*{-3pt}

\begin{lemma}\label{lem:34}
Let $(\lambda, f)$ be a solution to \eqref{P} and define $\Phi_0(z)$ and $\Phi_1(z)$ by  \eqref{Phi01}.
Then the functions $\Psi_0(z)$ and $\Psi_1(z)$, defined in \eqref{Psi01z}, satisfy \eqref{algebr} and solve the
integral equations \eqref{Psieqt}, where $\rho = \lambda^{1/(2\alpha)}$ and coefficients $c_0$ and  $c_1$ are
determined by \eqref{cf2c1} and
$
c_0 = \lim_{z\to \infty}\Psi_0(z).
$
\end{lemma} 

\subsection{Structure of the auxiliary system} 

Define the operator
\begin{equation}\label{A}
A f(t) := \frac 1 \pi \int_0^\infty \frac{e^{-\rho \tau}  }{\tau+t}\begin{pmatrix}
0 & g_0(\tau) \\
h_0(\tau) & 0
\end{pmatrix}
f(\tau) d\tau,
\end{equation}
which acts on functions $f:\Real_+\mapsto \Real^2$, and
consider the systems of integral equations
\vskip -10pt
\begin{equation}\label{pqr}
\begin{aligned}
p(t) =\, & A p(t) + \begin{pmatrix}
1 \\ 0
\end{pmatrix},
 \\
q(t) =\, & A q(t)   + \begin{pmatrix}
0 \\ 1
\end{pmatrix},
\\
r(t) =\, &
Ar(t) + \begin{pmatrix}
0 \\ t
\end{pmatrix},
\end{aligned}
\qquad t\in \Real_+.
\end{equation}
By changing variables, equations  \eqref{Psieqt} can be written as
$$
\begin{aligned}
\Psi_0(-\rho t) & =    \rho^\alpha\frac 1 \pi \int_0^\infty \frac{e^{-\rho \tau} g_0(\tau)}{\tau+  t}\Psi_1(-\rho\tau) d\tau + c_0,\\
\Psi_1(-\rho t) &   =  \rho^{-\alpha}\frac 1 \pi \int_0^\infty \frac{ e^{-\rho \tau} h_0(\tau) }{\tau+  t}\Psi_0(-\rho \tau)d\tau + c_1 (\rho b_\alpha-\rho t),
\end{aligned}
\qquad t\in  \Real_+,
$$
and therefore, by linearity,
\begin{equation}\label{eqPsi}
\begin{pmatrix}
\rho^{-\alpha}\Psi_0( \rho z)  \\
\Psi_1( \rho z)
\end{pmatrix}= \rho^{-\alpha} c_0 p(-z) + c_1 \rho b_\alpha q(-z) -c_1 \rho r(-z),
\end{equation}
where solutions to \eqref{pqr} are extended to $\mathbb C\setminus \Real_-$ by analyticity.
Plugging these expressions into \eqref{algebr} gives
\begin{equation}\label{cxieta}
c_0 \xi(\rho) + c_1 \eta(\rho) =0,
\end{equation}
where we denoted 
\begin{equation}\label{xieta}
\begin{aligned}
\xi(\rho) := &   X(\rho i)  p_1(-i) +\rho^{-\alpha}   e^{-\rho i}Y(-\rho i) p_2(i), \\
\nonumber
\eta(\rho) := &     X(\rho i)\rho^\alpha \Big(   \rho b_\alpha  q_1(-i) - \rho     r_1(-i)\Big)
+  e^{-\rho i}Y(-\rho i) \Big(  \rho b_\alpha   q_2(i) -  \rho   r_2(i) \Big).
\end{aligned}
\end{equation}
Equation \eqref{cxieta} has a nontrivial solution $(c_0,c_1)$ if and only if
\begin{equation}\label{cond_main}
\Im \big( \xi(\rho) \overline\eta(\rho)\big)=0.
\end{equation}
Thus we obtained the following converse to Lemma \ref{lem:34}.

\vspace*{-2pt}

\begin{lemma}\label{lem35}
Let $(p,q,r,\rho)$ be a solution to the system, which consists of integral equations
\eqref{pqr} and algebraic equation \eqref{cond_main}, such that $\rho >0$ and
\vskip -10pt
\begin{equation}\label{Apqr}
A\{p,q,r\}\in L_2(\Real_+;\Real^2).
\end{equation}
Define $\Psi_0(z)$ and $\Psi_1(z)$ by \eqref{eqPsi}, where $c_0$ is an arbitrary real constant and $c_1$ is determined
by \eqref{cxieta}. Let $\Phi_0(z)$ and $\Phi_1(z)$ be defined by equations \eqref{Psi01z} and $f$ be the inverse Laplace
transform of $\widehat f$, defined in \eqref{hatfrep} with $\lambda = \rho^{2\alpha}$ and $c_f$ as in \eqref{cf2c1}.
Then the pair $(\lambda, f)$ is a solution to \eqref{P}.
\end{lemma}

\subsection{Inversion of the Laplace transform} 

The inversion of Laplace transform, mentioned in Lemma \ref{lem35}, can be carried out by
integrating over the imaginary axis. In view of \eqref{hatfrep},
\begin{equation}\label{fx}
f(x)
=
-\frac 1 {2\pi i} \int_{-i\infty}^{i\infty}
\big(f_0(z) +f_1(z)\big)dz,
\end{equation}
\vskip -2pt \noindent 
where we defined
$$
f_0(z):= \frac { \Phi_0(z) } {\Lambda(z)} e^{zx} \quad \text{and}\quad f_1(z)
:=\Big(\frac { \Phi_1(-z) } {\Lambda(z)} +  \psi(1) \Big) e^{z(x-1)}.
$$
Integrating over suitable contours as in \cite[Lemma 5.8]{ChK} gives
\vskip -10pt
\begin{multline*}
\frac 1 {2\pi i} \int_{-i\infty}^{i\infty} \big(f_0(z) +f_1(z)\big)dz 
=   \Res(f_0, z_0) + \Res(f_0,-z_0) \\
+ \frac 1{2\pi i}\int_0^\infty \big(f_1^+(t)-f_1^-1(t)\big)dt
+ \frac 1 {2\pi i} \int_0^\infty \big(f_0^-(-t)-f_0^+(-t)\big)dt.
\end{multline*}
In view of symmetries \eqref{Lambdapm},
\vskip -11pt
\begin{align*}
f_0^-(-t)-f_0^+(-t)= &
-e^{-tx}\Phi_0(-t) \frac{2i \sin \big(\theta(t)-\alpha \pi\big)}{\big|\Lambda^+(t)\big|}, \\
f_1^+(t)-f_1^-(t)= &
-  e^{-t(1-x)} \Phi_1(-t)\frac {2i\sin \theta(t) } {|\Lambda^+(t)|}.
\end{align*}
The residues are complex conjugates,
$
\Res(f_0, -z_0) = \overline{\Res(f_0,z_0)}$,
and
$$
\Res(f_0,z_0) =
\frac { \Phi_0(z_0) } {\Lambda'(z_0)} e^{z_0x} = -\frac 1 {2i}
\rho^{1-\alpha} e^{-\frac \pi 2 \alpha i}\frac {\sin (\alpha \pi )} {\alpha \pi}
\Phi_0(z_0) e^{z_0x}.
$$

Plugging these expressions into \eqref{fx} we obtain
\begin{align}\label{feq}
&
f(x)  =\,    \Re\Big(\frac 1 { i}
\rho^{1-\alpha} e^{-\frac \pi 2 \alpha i}\frac {\sin (\alpha \pi )} {\alpha \pi}
\Phi_0(z_0) e^{z_0x}\Big) \\
&\nonumber
+ \frac 1{ \pi  }\int_0^\infty    e^{-t(1-x)} \Phi_1(-t)\frac { \sin \theta(t) } {|\Lambda^+(t)|}dt
+
\frac 1 { \pi  } \int_0^\infty  e^{-tx}\Phi_0(-t) \frac{  \sin \big(\theta(t)-\alpha \pi\big)}{\big|\Lambda^+(t)\big|} dt.
\end{align}

\subsection{Asymptotic analysis} 

Similarly to \cite[Lemma 5.6]{ChK}, the integral operator in \eqref{A} can be proved to be a contraction on
$L^2(\Real_+; \Real^2)$  for all sufficiently large values of $\rho$ and therefore equations \eqref{pqr}
have unique solutions with the property \eqref{Apqr}. As in \cite[Lemma 5.7]{ChK}, these solutions can be seen to
satisfy the estimates
\begin{equation}\label{pqr_est}
\left\|p(z) - \begin{pmatrix}
1 \\ 0
\end{pmatrix}\right\|\vee
\left\|q(z) - \begin{pmatrix}
0 \\ 1
\end{pmatrix}\right\|
\vee
\left\|r(z) - \begin{pmatrix}
0 \\ z
\end{pmatrix}\right\|
\le  C |z|^{-1} \rho^{-1}
\end{equation}
with a constant $C$.

\subsubsection{Eigenvalues} 
In view of \eqref{xieta} and estimates \eqref{pqr_est},
$$
\xi (\rho) \overline \eta(\rho) = (\rho i)^{-1}X_{c0}(i)
\rho^{\alpha}  (-i)^{\alpha-1}  e^{ \rho i}    X_{c0}( i) \big(    b_\alpha    +     i \big)\big(1+R(\rho)\big),
$$
where $R(\rho)$ can be shown to satisfy $|R(\rho)|\vee |R'(\rho)|\le C\rho^{-1}$ for some constant $C$.
Consequently  \eqref{cond_main} is equivalent to the equations
\begin{equation}\label{fla}
\begin{aligned}
\rho   =\,  & \pi n +\frac {\pi} 2 \alpha-   2  \arg\{X_{c0}(i)\}\\
&  -\arg\{    b_\alpha    +     i\}-  \mathrm{atan} \frac{\Im(R(\rho))}{1+\Re(R(\rho))}, \quad   n\in \mathbb Z.
\end{aligned}
\end{equation}

For each sufficiently large $n$ the derivative of the expression in the right hand side with respect to $\rho$
is less than unity. Hence if $n$ is taken large enough, so that, in addition, the operator \eqref{A} is a contraction,
the integro-algebraic system from Lemma \ref{lem35} has the unique solution,  obtained by
fixed-point iterations.
This, in turn, implies that the system \eqref{Psieqt}-\eqref{algebr} has the unique solution for all $\rho>0$
large enough, defined by \eqref{eqPsi}. Then since \eqref{P} has only finitely many eigenvalues on any bounded
interval, in view of Lemma \ref{lem:34}, all but possibly a finite number of eigenvalues correspond to solutions
of the system of Lemma \ref{lem35}.

In view of \eqref{fla}, the algebraic part of these solutions $\rho_n$  satisfies
$$
\rho_n  =\pi n +\frac {\pi} 2 \alpha-   2  \arg\{X_{c0}(i)\}-\arg\{    b_\alpha    +     i\} + O(n^{-1}).
$$
Here the last two quantities have explicit values, see Lemma \ref{lem-Xc0},
$$
\arg\big(X_{c0}(i)\big) =-\frac \pi 4 (1-\alpha)
\quad \text{and}\quad
\arg\big(b_\alpha+i\big) = \frac \pi {2\alpha},
$$
\vskip -3pt \noindent 
and consequently,
\vskip -12pt
\begin{equation}
\label{rhon}
 \rho_n
=\pi n  +  \frac \pi 2 (1 -1/\alpha)  + O(n^{-1}), \quad n\to\infty.
\end{equation}

Equations \eqref{fla} fix a particular enumeration of solutions to the integro-algebraic system
and, therefore, as argued above, of all sufficiently large eigenvalues $\lambda_n = \rho_n^{2\alpha}$.
Since any bounded interval contains only finitely many eigenvalues, the enumeration that puts all the eigenvalues
into the nondecreasing order may differ only by a finite shift.
This shift can be identified by comparing \eqref{rhon} at $\alpha=1$
to the classical problem \eqref{classics}, which shows that, in fact, the two enumerations coincide, that is,
\eqref{rho_n} holds. The technical details of this calibration procedure are the same as in \cite[Subsection 5.1.7]{ChK}.

\subsubsection{Eigenfunctions} 
When condition \eqref{cond_main} is satisfied, the equation \eqref{cxieta} determines the ratio
\vskip -12pt 
\begin{align*}
&
c_1/c_0 =   - \frac{\Re(\xi \overline \eta)}{|\eta|^2} =
- \Re\bigg(\frac{ \xi   }{ \eta }\bigg)  \\
&
= - \frac {\rho_n^{-1-\alpha}} {\sqrt{b_\alpha^2+1}}
\cos\Big(\rho_n  -\frac \pi 2 \alpha   +2\arg\{X_{c0}(i)\}-  \arg(b_\alpha    -     i)\Big)
 \big(1+O(\rho_n^{-1})\big) \\
&
=  - \frac {\rho_n^{-1-\alpha}} {\sqrt{b_\alpha^2+1}}
\cos (\pi n + O(n^{-1}))  \big(1+O(n^{-1})\big) \\
&
= - \frac {\rho_n^{-1-\alpha}} {\sqrt{b_\alpha^2+1}} (-1)^n
 \big(1+O(n^{-1})\big),
\end{align*}
\vskip -3pt \noindent
where we used estimates \eqref{pqr_est}.
Thus, in view of equation \eqref{eqPsi},
\vskip -12pt
$$
\Psi_0( i\rho_n  ) =    c_0 p_1(-i) + c_1    \rho_n^{1+\alpha} b_\alpha q_1(-i) -c_1   \rho_n^{1+\alpha}r_1(-i) =
c_0   \big(1+O(n^{-1})\big),
$$
\vskip -3pt \noindent
and  consequently,
\vskip -12pt 
\begin{align*}
\Phi_0(z_0) = \, & \Psi_0(i\rho_n) X(i\rho_n) = \Psi_0(i\rho_n) (i\rho_n)^{-1}X_{c,0}(i)  \\
&
=  c_0 \rho_n^{-1} e^{-\frac \pi 2 i}X_{c0}(i) \big(1+O(n^{-1})\big).
\end{align*}
The oscillating term in \eqref{feq} can be now simplified to
\vskip -12pt 
\begin{align*}
&
2\Re\Big(\frac 1 { 2i}
\rho^{1-\alpha} e^{-\frac \pi 2 \alpha i}\frac {\sin (\alpha \pi )} {\alpha \pi}
\Phi_0(z_0) e^{z_0x}\Big) \\
&
= \rho_n^{1-\alpha}\frac {\sin (\alpha \pi )} {\alpha \pi}  \Re\Big(
 e^{i \rho_n x-\frac \pi 2 (\alpha +1)i}\Phi_0(i\rho_n) \Big)  \\
&
= \rho_n^{ -\alpha}c_0|X_{c0}(i)|\frac {\sin (\alpha \pi )} {\alpha \pi}  \cos\Big(
   \rho_n x-\frac \pi 2 (\alpha +1) -\frac \pi 2  +\arg(X_{c0}(i))
\Big)\big(1+O(n^{-1})\big) \\
&
= \rho_n^{ -\alpha}c_0|X_{c0}(i)|\frac {\sin (\alpha \pi )} {\alpha \pi}  \cos\Big(
   \rho_n x- \frac \pi 4 \alpha -\frac {5\pi} 4
\Big)\big(1+O(n^{-1})\big)  \\
&
= -\rho_n^{ -\alpha}c_0\sqrt \alpha \frac {\sin (\alpha \pi )} {\alpha \pi}  \sin\Big(
   \rho_n x  +\frac { \pi} 4 (1-\alpha) \Big)\big(1+O(n^{-1})\big),
\end{align*}
\vskip-3pt \noindent 
where we used \eqref{Xc0val}. The boundary layer terms are treated similarly,
\vskip -12pt
\begin{align*}
&
\frac 1{ \pi  }\int_0^\infty    e^{-t(1-x)} \Phi_1(-t)\frac { \sin \theta(t) } {|\Lambda^+(t)|}dt  \\
&
= \rho_n^{1-\alpha}\frac {\sin (\alpha \pi )} {\alpha \pi}\frac \alpha { \pi  }\int_0^\infty  \Phi_1(-\tau \rho_n)\frac { \sin \theta_0(\tau  ) } {\gamma_0(\tau)} e^{-\rho_n  (1-x)\tau}d\tau  \\
&
=  (-1)^n\frac {c_0} {\rho_n^\alpha}   \frac {\sin (\alpha \pi )} {\alpha \pi}\frac \alpha { \pi  } \! \int_0^\infty
 \frac {\tau^\alpha(\tau \!-\! b_\alpha) } {\sqrt{b_\alpha^2\!+\!1}}\frac{ X_{c0}(-\tau)}{\tau}
\frac { \sin \theta_0(\tau  ) }
{\gamma_0(\tau)} e^{-\rho_n  (1-x)\tau}d\tau\big(1\!+\!O(n^{-1})\big),
\end{align*}
\vskip -4pt \noindent 
where $\gamma_0(\tau)$ is defined by the equation
\vskip -12pt
\begin{equation}\label{gamma0}
|\Lambda^+(\tau \rho_n)| = \rho_n^\alpha
\frac \pi {\sin (\alpha \pi )} \Big|
 \tau^{\alpha } -      e^{ \alpha \pi  i}\tau^{-\alpha}
 \Big| =: \rho_n^\alpha \frac \pi {\sin (\alpha \pi )}\gamma_0(\tau),
\end{equation}
\vskip-4pt \noindent 
and we used the uniform approximation due to \eqref{pqr_est},
\vskip-12pt 
\begin{align*}
&
\Phi_1(-\tau \rho_n) =  Y(-\tau \rho_n)\Psi_1(-\tau \rho_n)  \\
&
= (\tau \rho_n)^{\alpha-1} X_c(-\tau \rho_n)
\Big(
\rho_n^{-\alpha} c_0 p_2(\tau) + c_1 \rho_n b_\alpha q_2(\tau) -c_1 \rho_n r_2(\tau)
\Big) \\
&
= - c_0\rho_n^{ -1} (-1)^n \frac 1 {\sqrt{b_\alpha^2+1}}\tau^{\alpha-1} X_{c0}(-\tau)
   (b_\alpha  -  \tau)
  \big(1+O(n^{-1})\big).
\end{align*}

Similarly,
\begin{align*}
&
\frac 1 { \pi  } \int_0^\infty  e^{-tx}\Phi_0(-t) \frac{  \sin \big(\theta(t)-\alpha \pi\big)}{\big|\Lambda^+(t)\big|} dt\\
&
= \rho_n\frac 1 { \pi  } \int_0^\infty  e^{-\rho_n  x \tau}\Phi_0(-\rho_n \tau ) \frac{  \sin \big(\theta(\rho_n \tau)-\alpha \pi\big)}{\big|\Lambda^+(\rho_n \tau)\big|} d\tau \\
&
= \rho_n^{1-\alpha}\frac  {\sin (\alpha \pi )}{\pi\alpha} \frac \alpha { \pi  } \int_0^\infty  e^{-\rho_n  x \tau}\Phi_0(-\rho_n \tau ) \frac{  \sin \big(\theta_0( \tau)-\alpha \pi\big)}{ \gamma_0(\tau)} d\tau \\
&
= - c_0 \rho_n^{ -\alpha}\frac  {\sin (\alpha \pi )}{\pi\alpha} \frac \alpha { \pi  } \int_0^\infty
     \frac{X_{c0}(-\tau)}{\tau}   \frac{  \sin \big(\theta_0( \tau)-\alpha \pi\big)}{ \gamma_0(\tau)}  e^{-\rho_n  x \tau}d\tau(1+O(n^{-1})),
\end{align*}
where we used the approximation
\begin{align*}
\Phi_0(-\rho_n \tau )= & X(-\rho_n \tau)\Psi_0(-\rho_n \tau)  \\
&
= X(-\rho_n \tau) \Big(
c_0 p_1(\tau) + c_1 \rho_n^{1+\alpha} b_\alpha q_1(\tau) -c_1 \rho_n^{1+\alpha} r_1(\tau)
\Big) \\
&
=  - c_0 \rho_n^{-1}  \tau^{-1}X_{c0}(-\tau)   (1+O(n^{-1})).
\end{align*}
Assembling all parts in \eqref{feq} together and normalising up to the unit norm, we arrive at
\eqref{fn}
with
\begin{equation}\label{Upsilons}
\begin{aligned}
\Upsilon_0(\tau) & =
\frac {\sqrt {2\alpha}} { \pi  } \frac{X_{c0}(-\tau)}{\tau}   \frac{  \sin \big(\theta_0( \tau)-\alpha \pi\big)}{ \gamma_0(\tau)},\\
\Upsilon_1(\tau) & =
\frac {\sqrt {2\alpha}} { \pi  }\frac {\tau^{\alpha }(b_\alpha  -  \tau)} {\sqrt{b_\alpha^2+1}} \frac{X_{c0}(-\tau)}{\tau}
\frac { \sin \theta_0(\tau  ) }
{\gamma_0(\tau)},
\end{aligned}
\end{equation}
where, cf. \eqref{gamma0} and  Lemma \ref{lem-Lambda} (c),
$$
\gamma_0(\tau)  = \Big(\tau^{2\alpha} -2\cos(\alpha\pi) +\tau^{-2\alpha}\Big)^{1/2}
 \quad \text{and}\quad
\theta_0(\tau) =  -\mathrm{atan}\frac{\sin (\alpha \pi)}{\tau^{2\alpha}-\cos(\alpha \pi)},
$$
and, cf. Lemma \ref{lem-Xc0}, $b_\alpha = \cot  \left(\frac \pi {2\alpha}\right)$ and
$$
X_{c0}(-\tau)= \exp \left(\frac 1 \pi \int_0^\infty \frac{\theta_0(t)}{t+\tau}dt\right).
$$

\section{Conclusions} 

\setcounter{section}{4}
\setcounter{equation}{0}\setcounter{theorem}{0}

This paper constructs an asymptotic approximation of solutions to a fractional Sturm-Liouville problem of a particular type.
The obtained eigenvalues asymptotics \eqref{rho_n} is exact up to the second order, with an estimate for the residual.
This precision turns out to be sufficient for uniform approximation \eqref{fn} of the eigenfunctions, which sheds some
light on their general structure. The suggested method is based on reduction to a Hilbert boundary value problem.
While the details are worked out in a specific problem, the technique is expected to be applicable
to other fractional eigenproblems of interest.

\vspace*{-4pt}

\section*{Acknowledgements} 

P. Chigansky’s research was funded by ISF 1383/18 grant.




 \bigskip  

 \it

 \noindent
$^1$ Department of Statistics, 
The Hebrew University of Jerusalem \\
Mount Scopus, Jerusalem 91905, 
ISRAEL  \\[3pt]
  e-mail: pchiga@mscc.huji.ac.il \ (Corresponding author)\\
  \hspace*{1cm} 
\hfill Received: June 12, 2020,\ Revised: ......, 2021  \\[6pt] 
$^2$ Laboratoire de Statistique et Processus\\
Le Mans Universit\'{e}, FRANCE \\[3pt]
  e-mail: marina.kleptsyna@univ-lemans.fr

\end{document}

%% file: fcaa_style_DG.tex
\usepackage{graphicx}
\usepackage{epsfig}

\usepackage{amsthm}
\usepackage{amsmath}
\usepackage{latexsym}
\usepackage{amsfonts}
\usepackage{amssymb}

\newtheoremstyle{theorem}
  {15pt}          
  {15pt}  
  {\sl}  
  {\parindent}
  {\sc}  
  {. }   
  { }    
  {}     
\theoremstyle{theorem}
\newtheorem{lemma}{Lemma}[section]
\newtheorem{theorem}{Theorem}[section]

\newtheoremstyle{defi}
  {15pt}          
  {15pt}  
  {\rm}  
  {\parindent}     
  {\sc}  
  {. }    
  { }    
  {}     
\theoremstyle{defi}

 \def\proofname{\indent\rm P\ r\ o\ o\ f.}
 \def\proofend{\hfill$\Box$}
 \def\qed{\hfill$\Box$} 